\documentclass{amsart}
\usepackage{amsmath}
\usepackage{amsfonts}
\usepackage{amssymb}
\usepackage{graphicx}
\usepackage{multirow}
\usepackage{stix}
\usepackage{courier}
\usepackage{wrapfig}
\usepackage{eucal}
\usepackage{tikz}
\newcommand{\mysetminus}{\hbox{\tikz{\draw[line width=0.6pt,line cap=round] (3pt,0) -- (0,6pt);}}}

\numberwithin{equation}{section}
\begin{document}

\title[Quadri-Figures in Cayley-Klein Planes II: The Miquel-Steiner Theorem]{Quadri-Figures in Cayley-Klein Planes II:\\The Miquel-Steiner Theorem}
\author{Manfred Evers}
\curraddr[Manfred Evers]{Bendenkamp 21, 40880 Ratingen, Germany}
\email[Manfred Evers]{manfred\_evers@yahoo.com}
\date{\today}

\begin{abstract}
The Miquel-Steiner theorem for a quadrilateral in the euclidean plane states that the circumcircles of the four component triangles intersect at a single point, which now is called the \textit{Miquel-Steiner point} of the quadrilateral.\\ 
The Miquel-Steiner theorem for euclidean planes needs to be slightly modified for Minkow\-ski and galilean planes: Either the circumcircles of the four component triangles touch each other at an isotropic point, or they intersect transversally at an anisotropic point.\\
In elliptic and hyperbolic planes, as well as in dual-euclidean and dual Minkow\-ski planes, the Miquel-Steiner theorem does not hold in this form. Instead, a weaker version applies: The circumcircles of the four component triangles of a quadrilateral have a common radical center, which we will also call the Miquel-Steiner point.\\
For specific quadrilaterals (such as cyclic quadrilaterals), the location of the Miquel-Steiner point can be determined more precisely.
\end{abstract}
\maketitle \hspace*{\fill}\vspace*{-3 mm}\\
\section*{Introduction}

In order to introduce the terminology and fix notations, some definitions, rules and basic theorems are given in the first section.\\
\hspace*{2 mm}In the second section, we look at how far the Miquel-Steiner theorem, known from euclidean geometry, and theorems based on it can be transferred to elliptic, hyperbolic, dual-euclidean and dual Minkow\-ski  geometry.\\
\hspace*{2 mm}In the third section, we give an analytical proof of the Miquel-Steiner theorem in metric-affine planes (euclidean, Minkow\-ski and galilean planes). This may seem surprising at first, since there are already short and elegant synthetic proofs for euclidean geometry. However, these are mostly based on the inscribed angle theorem. We show that a transfer to Minkow\-ski geometry is possible.\vspace*{0.5 mm}\\
\hspace*{2 mm}The first three sections treat circles as regular conics, while the fourth briefly extends this to singular circles. We provide examples showing that the Miquel-Steiner theorem holds also for these cases, though a general proof is omitted. \vspace*{0.5 mm}\\
\section{Fundamentals}

In the following, we consider the real vector space $\mathbf{U}:=\mathbb{R}^3$ embedded in the complex vector space $\mathbf{V}:=\mathbb{C}^3$ and the real projective plane $\mathrm{P}\mathbf{U}$  embedded in the complex projective plane $\mathrm{P}\mathbf{V}$.
The elements of $\mathrm{P}\mathbf{V}$ are called points; and if $P = \mathbb{C}\mathbf{v}$ is such a point, we write $P=[\mathbf{v}]$. Points in \(\mathrm{P}\mathbf{U}\) are called \textit{real}.\\
\hspace*{2 mm}Lines  in $\mathrm{P}\mathbf{V}$ are zero sets of projective linear forms; for any line $\mathcal{L}$, there exists a linear form $l{:\,}\mathbf{V} \to \mathbb{C}$ such that $\mathcal{L}=\{[\mathbf{v}]\in  \mathrm{P}\mathbf{V}\,|\,l(\mathbf{v})=0\}$. The line passing through two distinct points $P_1$ and $P_2$ is denoted by $P_1\times P_2$, the intersection point of two lines $\mathcal{L}_1$ and $\mathcal{L}_2$ by $\mathcal{L}_1{\;\times\;}\mathcal{L}_2$. A line $\mathcal{L} \subset \mathrm{P}\mathbf{U}$ is called a \textit{real line}.\\
\hspace*{2 mm}A conic $\mathcal{K}\subset \mathrm{P}\mathbf{V}$ is the zero set of a non-zero projective quadratic form $q{:\,}\mathbf{V} \to \mathbb{C}$. The conic $\mathcal{K}$ defines the quadratic form $q$ up to a non-zero complex factor. A projective quadratic form $q{:\,}\mathbf{V} \to \mathbb{C}$ is called \textit{real} if it maps elements in $\mathbf{U}$ to real numbers. A conic is defined to be real if it is the zero set of a non-zero real projective quadratic form and contains at least two real points.\\

\noindent\textit{Semi Cayley-Klein planes}, cf. \cite{Ev2}.\\
We now concentrate on the real projective plane $\mathrm{P}\mathbf{U}$ and on real quadratic forms $\phi{:} \mathbf{V} \to \mathbb{C}$. We select such a form $\phi_{\circ}$, call it the \textit{absolute form} and call the conic $\mathcal{K}_{\circ}:= \{\,[\mathbf{p}]\,|$ $\,\phi_{\circ}(\mathbf{p})=0\,\}$ the \textit{absolute conic}.
$\mathrm{P}\mathbf{U}$  becomes a \textit{semi-Cayley-Klein plane} by selecting from the set of automorphisms on $\mathrm{P}\mathbf{V}$
only those that leave $\mathrm{P}\mathbf{U}$ and the absolute conic invariant. We call these maps \textit{semi-isometries}.\\ 
Based on the signature of the quadratic form $\phi_{\circ}$, we can make the following classification of the geometries on $\mathrm{P}\mathbf{U}$. The geometry of the plane is elliptic, hyperbolic, dual-euclidean, dual Minkow\-ski, or metric-affine, depending on whether this signature is $(3,\!0,\!0),(2,\!1,\!0),$ $(1,\!1,\!0),$ $(1,-1,0)$, or $(1,0,0)$.\\
For elliptic and hyperbolic planes, semi-isometries are isometries. This is not true for the non-regular planes.\\

\noindent\textit{Congruence of anisotropic points}\\
We define two anisotropic points $P,Q \in \mathrm{P}\mathbf{U}$ to be congruent (we write $P\cong  Q$) if there is an semi-automorphism on $\mathrm{P}\mathbf{V}$ that maps $P$ to $Q$. 
It can be shown that two anisotropic points are congruent if and only if they are in the same connected component of $\mathrm{P}\mathbf{U}\,\mysetminus\,\mathcal{K}_{\circ}$.\\
In elliptic, metric-affine and dual-euclidean planes, any two anisotropic points are congruent. In hyperbolic and dual Minkow\-ski planes, however, there are two congruence classes.\\

\noindent\textit{Line segments and their midpoints.}\\
We assign a vector $P^\circ\in \boldsymbol{U}$  to each anisotropic point $P \in \mathrm{P}\mathbf{U}$. \\
First, we define a function $\chi{:}\;\mathbb{R}^3 \rightarrow \{-1,0,1\}$ by\\ \vspace*{2.5 mm}  
{\!$\chi(p_1,p_2,p_{3}) = 
\begin{cases}
 \;\,0,\;\text{if } (p_1,p_2,p_3) = (0,0,0)\;, \\
 \;\,1,\;\text{if } (p_1,p_2,p_3) > (0,0,0)\; \textrm{with respect to the lexicographic order,} \\
-1,\text{if } (p_1,p_2,p_3) < (0,0,0)\; \textrm{with respect to the lexicographic order}; \\
\end{cases}$}\vspace*{-1.5 mm}\\
then we set $\displaystyle P^\circ := \frac{\chi(\boldsymbol{p})}{\sqrt{|\phi_{\circ}(\boldsymbol{p})|}} \boldsymbol{p}$.\vspace*{+1.0 mm}\\
\hspace*{2 mm}For real points $P,Q,R$ and complex numbers $\alpha,\beta,\gamma$, we write $\alpha P{+}\beta\,Q{+}\gamma R$ 
more concisely instead of $[\alpha P^\circ{+}\beta\,Q^\circ{+}\gamma R^\circ]$.\vspace*{3 mm}\\
\hspace*{2 mm}The two \textit{line segments} bounded by anisotropic real points $P$ and $Q$ are\vspace*{0.5 mm}\\ 
\centerline{$\big[P,Q\big]_+{\,:=\,}$$\{\alpha P{+}\beta\,Q\, |\,\alpha,\beta \in \mathbb{R}_+\}$ and $\big[P,Q\big]_-{\,:=\,}\{\alpha P{+}\beta\,Q\,|\, \alpha,\beta \in \mathbb{R}\,,\,\alpha\beta\leq 0\}.$} \vspace*{0.5 mm}
These two segments have midpoints precisely when $P\cong Q$.\\
If $P\cong Q$, the midpoint of $\big[P,Q\big]_\pm$ is $P\pm Q$. At least one of these two midpoints is anisotropic.\\

\noindent\textit{Triangles.}\\ Given three non-collinear anisotropic real points $R,S,T$, there are four distinct triangles with vertices $R,S,T$,\\
\hspace*{10 mm}$\Delta_0(R,S,T){\,{{\scriptstyle{:}}\!\!=}\,}\{(q_1R+q_2S+q_3T)\in \mathrm{P}\mathbf{U}\,|\,q_1,q_2,q_3\geq 0\}$,\\
\hspace*{10 mm}$\Delta_1(R,S,T)$ ${\,{{\scriptstyle{:}}\!\!=}\,}\{(q_1R+q_2S+q_3T)\in \mathrm{P}\mathbf{U}\,|\,q_1\geq 0,{q_2,q_3\leq 0}\}$, $\,\dots\,$.\vspace*{1 mm}\\
In the following, "triangle $RST$" always refers to "$\Delta_0(R,S,T)$".\vspace*{2 mm}\\
\noindent\textit{Barycentric coordinates with reference to a triangle $ABC$.}\\
We only want to consider triangles that have a circumcircle. This is only the case if all three (anisotropic) vertices belong to the same congruence class. We select such a triangle $ABC$; it will be our reference triangle.\\ 
We calculate its centroid:\\
$\hspace*{6 mm}G = ((A+B)\times C)\times ((B+C)\times A) = [(A^\circ{\,+\,}B^\circ)\times C^\circ)\times ((B^\circ{\,+\,}C^\circ)\times A^\circ)]\\
\hspace*{8.5 mm}= [A^\circ + B^\circ + C^\circ] = A+B+C.$\\
Given any point $P=[\boldsymbol{p}]\in \mathrm{P}\mathbf{V}$, we write $P = [p_A{:}p_B{:}p_C]$ if there exists a number $t \in \mathbb{C}$ such that $\boldsymbol{p} = t\,(p_A A^{\circ} + p_B B^{\circ} + p_C C^{\circ})$. $[p_A{:}p_B{:}p_C]$ is a representation of $P$ in \textit{homogeneous barycentric coordinates} with respect to $(A,B,C)$.\vspace*{1 mm}\\
Here are some examples:\vspace*{0.5 mm}\\
\hspace*{6 mm}$A =[1{:}0{:}0]\,,\;B+C = [0{:}1{:}1]\,,\;G=[1{:}1{:}1]\,.$\vspace*{0.5 mm}\\
\hspace*{2 mm}The absolute conic is\\
\hspace*{6 mm}$\mathcal{K}_{\circ}=\{ [p_A{:}p_B{:}p_C] \in \mathrm{P}\mathbf{V}\;|\,p_A^2+p_B^2+p_C^2+2r_Ap_Bp_C+2r_Bp_Cp_A+2r_Cp_Ap_B=0\}$, \\
\hspace*{2 mm}with suitable real numbers $r_A,r_B,r_C$. \vspace*{2 mm}\\
\textit{Circles}\\Circles are special nonsingular real conics. We distinguish two cases.\\
\hspace*{2 mm}In metric affine planes, the absolute conic is the double line $\{ [p_A{:}p_B{:}p_C] \in \mathrm{P}\mathbf{V}\;|\,(p_A+p_B+p_C)^2=0\}$. The line at infinity, $\mathcal{L}_{\infty} := \{ [p_A{:}p_B{:}p_C] \in \mathrm{P}\mathbf{V}\;|\,p_A+p_B+p_C=0\}$, is the tripolar of the centroid $G$ of triangle $ABC$. To determine when a non-singular real conic is a circle, we define a semi-Cayley-Klein structure on the line at infinity by speci\-fying two (not necessarily distinct) isotropic points $I,J$ on $\mathcal{L}_{\infty}$. However, neither of these two points may lie on a side line of triangle ABC. A nonsingular real conic is a circle if and only if it shares exactly these points, the so-called circular points of the plane, with the line at infinity.
We distinguish the following cases: If there are two real circular points, the geometry of the plane is a Minkow\-ski geometry. If there are circular points that lie in $\mathrm{P}\mathbf{V}\,\mysetminus\,\mathrm{P}\mathbf{U}$, the geometry of the plane is euclidean. In a galilean plane, there is just one circular point; this lies in $\mathrm{P}\mathbf{U}$.\vspace*{0.5 mm}\\
\hspace*{2 mm}In non-affine Cayley-Klein planes, a nonsingular real conic $\mathcal{C}$ that differs from the absolute conic is a \textit{circle} if there exists a line $\mathcal{L}$ in $\mathrm{P}\mathbf{U}$ whose anisotropic points are all symmetry points of $\mathcal{C}$. A nonsingular real conic, other than a circle, can have at most three real symmetry points. The pole of $\mathcal{L}$ with respect to $\mathcal{C}$ is called the circle center.  If this circle center is an anisotropic point, then it is also a symmetry point of the circle.\\
In dual-euclidean and in dual Minkow\-ski planes, all circles have the same center O. It is the \textit{absolute point at infinity}, the singular point on the absolute conic.\\

\noindent\textit{The circumcircle of triangle} $ABC$.\\
\hspace*{2 mm}We begin by considering metric affine planes. A nonsingular real circumconic of triangle \(ABC\) can be defined by the equation\vspace*{0.5 mm}\\
\centerline{\(k_{A}z_{B}z_{C}+k_{B}z_{C}z_{A}+k_{C}z_{A}z_{B}=0\)\,,}\vspace*{0.5 mm}
where \(k_A, k_B, k_C \in \mathbb{R}\,\mysetminus\,\{0\}\). The point \([k_A{:}k_B{:}k_C]\) is the perspector of this circumconic.\\
The circumcircle of \(ABC\), denoted by \(\mathcal{C}_{ABC}\), is the unique circumconic that passes through the circular points at infinity. The perspector of \(\mathcal{C}_{ABC}\) is known as the Lemoine point of the triangle. The center of the circumcircle is given by the cevian quotient of the Lemoine point \(K\) and the centroid \(G\), yielding the point\vspace*{0.5 mm}\\
\centerline{\(O= [k_{A}(-k_{A}+k_{B}+k_{C}):k_{B}(k_{A}-k_{B}+k_{C}):k_{C}(k_{A}+k_{B}-k_{C})]\)\,.}
Consequently, the circular points \(I\) and \(J\) are the points \(P=[z_A{:}z_B{:}z_C]\) in \(\mathrm{P}\mathbf{V}\) that simultaneously satisfy the linear and quadratic equations\vspace*{0.5 mm}\\
\centerline{\(z_{A}+z_{B}+z_{C}=0\quad \text{and}\quad k_{A}z_{B}z_{C}+k_{B}z_{C}z_{A}+k_{C}z_{A}z_{B}=0\)\,.}\vspace*{2 mm}\\
\hspace*{2 mm}We now turn to non-affine planes. The circumcircle of triangle \(ABC\) is the conic\\ 
\centerline{$\hspace*{2 mm}\mathcal{C}_{ABC}\,=\,\{\,[p_A{{:}}p_B{{:}}p_C] \in \mathrm{P}\mathbf{V}\,|\,(r_A-1)p_Bp_C+(r_B-1)p_Cp_A+(r_C-1)p_Ap_B=0\}$\,,} 
where $r_A,r_A,r_B$ are the real numbers introduced in the penultimate subsection. The perspector of this circle is the Lemoine point $K{\,=\,}[k_A{:}k_B{:}k_C]{\,=\,}[r_A-1{\,:\,}r_B-1{\,:\,}r_C-1]$, and the point $O=[o_A{:}o_B{:}o_C]=[k_A(k_A-k_B-k_C):k_B(k_B-k_C-k_A):k_C(k_C-k_A-k_B)]$ is the \textit{circumcenter} of triangle $ABC$. The polar of $O$ with respect to the circumcircle is identical to the tripolar of $G$  with respect to $ABC$. \\
In elliptic and hyperbolic planes, this line is also the polar of $O$ with respect to $\mathcal{K}_{\circ}$. \\

\noindent\textit{Congruence of lines.}\\
In hyperbolic and Minkowski planes, anisotropic real lines fall into two congruence classes; in all other cases, they are all congruent. In a hyperbolic plane, two such lines are congruent if and only if their poles with respect $\mathcal{K}_{\circ}$ are congruent. In Minkowski planes, $\mathcal{L}_\infty$ and all real lines through a circular point are isotropic; and two anisotropic lines are congruent if and only if their intersections with $\mathcal{L}_{\infty }$ lie in the same connected component of $\mathcal{L}_\infty\mysetminus{\,\mathcal{C}_{ABC}}$. \\

\noindent\textit{Orthogonality of points and lines.}\\
\hspace*{2 mm}First, we consider the case of a regular absolute conic $\mathcal{K}_{\circ}$. 
If  $P$ is a point and \(\mathcal{L}\) is a line in \(\mathrm{P}\mathbf{U}\), we denote the polar line of $P$ with respect to $\mathcal{K}_{\circ}$ by $P^\perp$ and the pole of \(\mathcal{L}\) by \(\mathcal{L}^\perp\). Now we define the orthogonality of a pair $(P,\tilde{P})$ of points and a pair $(\mathcal{L},\tilde{\mathcal{L}})$ of lines by\\
\centerline{\(P \perp \tilde{P}\,:\Leftrightarrow P\,\in \tilde{P}^\perp\;\;\) and \(\;\;\mathcal{L} \perp \tilde{\mathcal{L}}\,:\Leftrightarrow \,\mathcal{L}^\perp \in \tilde{\mathcal{L}}\;\;.\;\;\) }
\hspace*{2 mm}In euclidean and Minkowski planes, the orthogonality of a pair $(\mathcal{L},\tilde{\mathcal{L}})$ of anisotropic lines is defined by: 
{\hspace*{18 mm}$\mathcal{L}\perp\tilde{\mathcal{L}} \,:\Leftrightarrow\;\operatorname{cr}(\mathcal{L},\tilde{\mathcal{L}};I,J)=-1\;\;$}\\    ($\operatorname{cr}$ stands for \textit{cross ratio}, and $I,J$ are the circular points).  \\
\hspace*{2 mm}In dual-euclidean and dual-Minkowski planes we define the orthogonality of a pair $({P},\tilde{P})$ of anisotropic points by: 
$\hspace*{2 mm}P\perp\tilde{P} \,:\Leftrightarrow \;\operatorname{cr}({M_\infty\times P},{M_\infty\times \tilde{P}};\mathcal{I},\mathcal{J})=-1$, \\
where $\mathcal{I},\mathcal{J}$ are the two lines with $\mathcal{I}\cup \mathcal{J}=\mathcal{K}_{\circ}$ and $M_\infty=\mathcal{I}\times \mathcal{J}$ is the (real) point at infinity. \\
\hspace*{2 mm}In all these cases, orthogonality is a symmetric relation.\\

\noindent\textit{Quadrilaterals}\\ 
A quadrilateral in the projective plane $\mathrm{P}\mathbf{U}$ consists of four lines in a general position, no three of which are concurrent.
Any three of these lines are the side-lines of a triangle, which is called a component-triangle of the quadrilateral. The quadrilateral
together with its four component triangles is called a \textit{complete quadrilateral}. 
For a  complete quadrilateral in a Cayley-Klein plane $\mathrm{P}\mathbf{U}$ we also require that the vertices of all component triangles are anisotropic points that all belong to the same congruence class, so every component triangle has its own circumcircle.\vspace*{-2 mm}

\section{The Miquel-Steiner theorem for non-affine Cayley-Klein planes.}

Throughout the entire section, $ABC$ is the reference triangle, and the absolute conic is 
$\mathcal{K}_{\circ}=\{ [p_A{:}p_B{:}p_C] \in \mathrm{P}\mathbf{V}\;|\,p_A^2+p_B^2+p_C^2+2r_Ap_Bp_C+2r_Bp_Cp_A+2r_Cp_Ap_B=0\}$.\vspace*{-2 mm}\\
\begin{figure}[!hbpt]
\includegraphics[height=7.3cm]{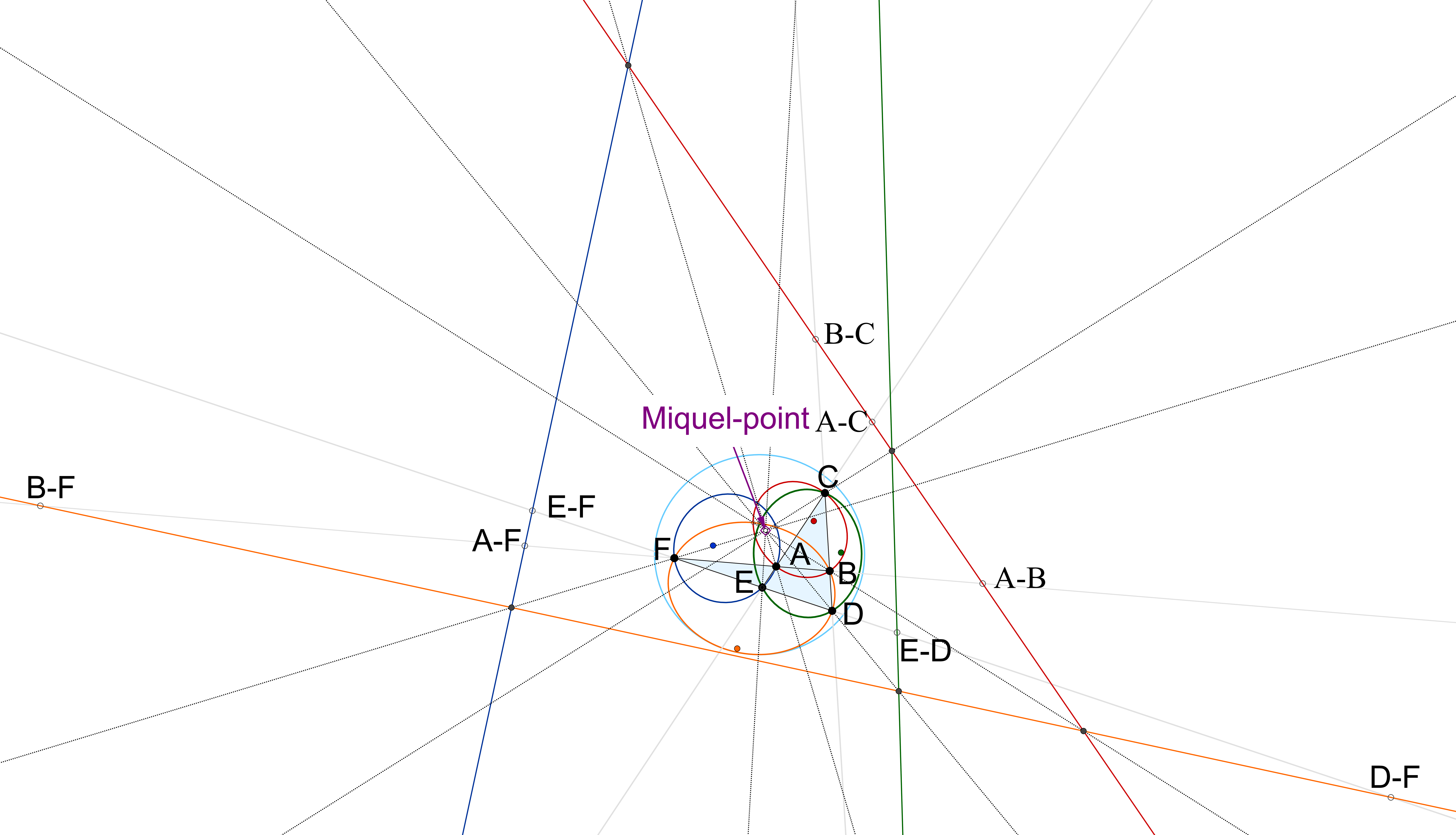}
\centering
\caption{The Miquel-Steiner point of a complete quadrilateral in a hyperbolic plane. The light blue circle is the absolute conic in this plane.}
\vspace*{2 mm}\end{figure}

\begin{figure}[!htbp]
\includegraphics[height=7.1cm]{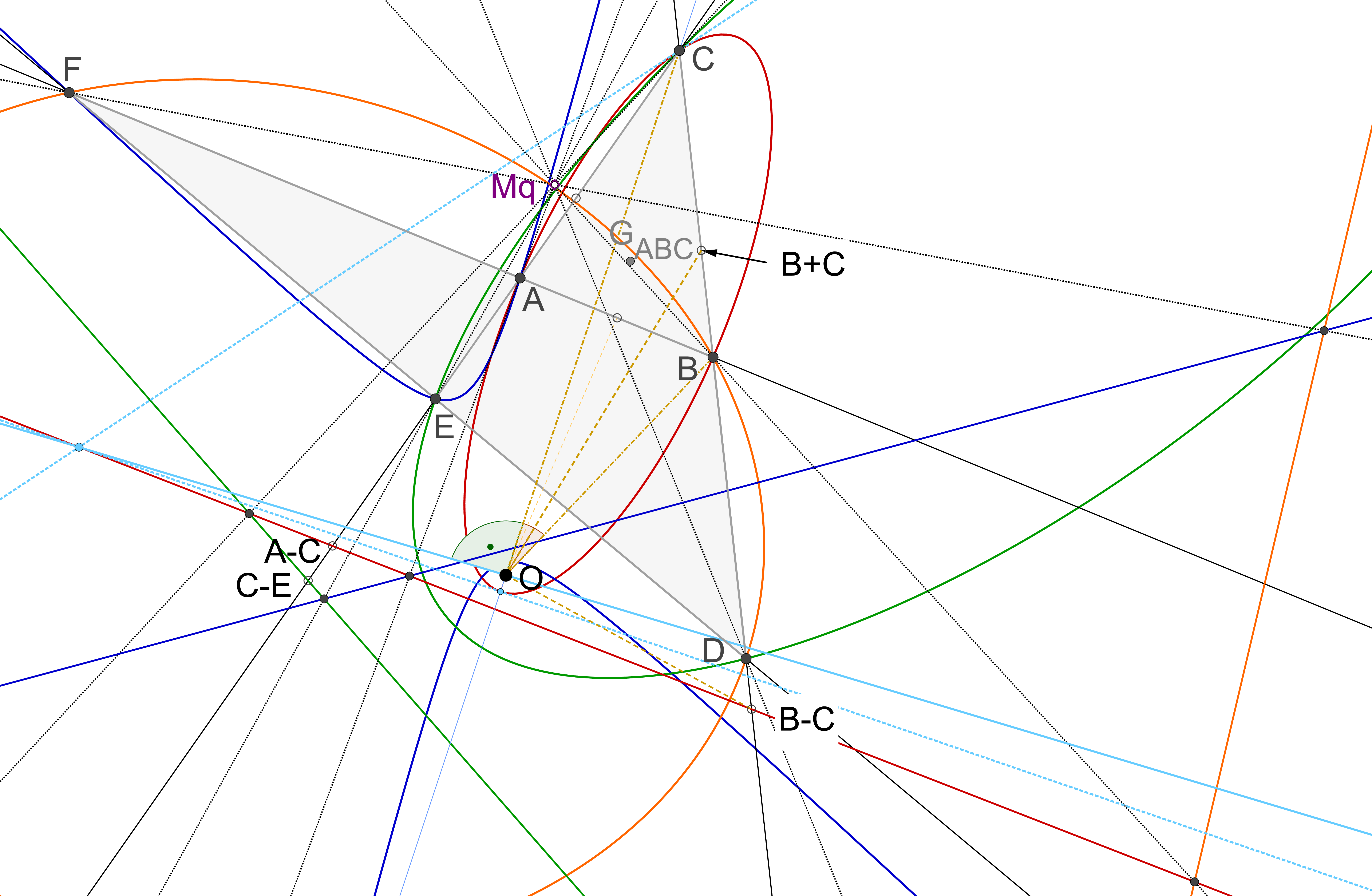}
\centering
\caption{The Miquel-Steiner point of a complete quadrilateral in a dual-euclidean plane. The isotropic point $O$ is the absolute point at infinity.}
\vspace*{2 mm}\end{figure} \hspace*{-1 mm}

\noindent\textit{Theorem} 1.\;\;
The circumcircles of the four component triangles of a complete quadrilateral have a common radical center, which we call the \textit{Miquel-Steiner point} of the quadrilateral. \vspace*{1 mm} \\
\noindent{\textit{Proof}.} 
The vertices of the complete quadrilateral are labeled as shown in Figure 1.
For the points $D, E, F$ we can find positive real numbers $l,m,n$ such that $D=[0{:}n{:}{-}m], E=[n{:}0{:}{-}l], F=[m{:}{-}l{:}0]$,\\
This results in the vectors $A^\circ,B^\circ,C^\circ,D^\circ,E^\circ,F^\circ$:\vspace*{1 mm} \\
\hspace*{10 mm}$A^\circ=(1,0,0),B^\circ=(0,1,0),C^\circ=(0,0,1),\vspace*{1 mm}\\
\hspace*{10 mm}D^\circ=(0,dn,{-}dm),E^\circ=(en,0,{-}el),F^\circ=(fm,{-}fl,0)$ \vspace*{1 mm}\\
\hspace*{2 mm}with $d{\,:=\,}1/\sqrt{m^2{+}n^2{-}2r_Amn\,}, e{\,:=\,}1/\sqrt{l^2{+}n^2{-}2r_Bln\,}, f{:=\,}1/\sqrt{l^2{+}m^2{-}2r_Clm\,}$.\vspace*{1 mm}\\
Each of the four component triangles of the quadrilateral has its own circumcircle, and any three of these circles have there own radical center.\vspace*{0.5 mm}\\
Let $\mathcal{C}_1, \mathcal{C}_2, \mathcal{C}_3, \mathcal{C}_4$ be the circumcircles of the component triangles $ABC, AEF, CDE$ and $BDF$, with circumcenters $O_1, O_2, O_3, O_4$, respectively. 
To find the radical center of $\mathcal{C}_{1}, \mathcal{C}_{2}$ and $\mathcal{C}_{3}$, we determine the intersection of their three radical axes:\\
\centerline{$\mathcal{L}_A := (O_1\times O_2)^\perp\times A = \big(\big((A{-}B)\times(A{-}C)\big)\times \big((A{-}E)\times(A{-}F)\big)\big)\times A$,}
\centerline{$\mathcal{L}_B := (O_1\times O_3)^\perp\times B = \big(\big((A{-}B)\times(A{-}C)\big)\times \big((B{-}D)\times(B{-}F)\big)\big)\times B$,}
\centerline{$\mathcal{L}_C := (O_1\times O_4)^\perp\times C = \big(\big((A{-}B)\times(A{-}C)\big)\times \big((C{-}D)\times(C{-}E)\big)\big)\times C$.}
Computing the intersection yields the radical center $M\!q$:\\
\centerline{$M\!q =\mathcal{L}_A \times \mathcal{L}_B =  \mathcal{L}_A \times \mathcal{L}_C =[ef(d(m{-}n){\,-\,}1):df(1{\,-\,}e\,(l{-}n)):de\,(f (l{-}m){\,-\,}1)].$}
It can be easily checked that this point $M\!q$ is also a point on the lines $\mathcal{L}_D := (O_3\times O_4)^\perp\times D$, $\mathcal{L}_E := (O_2\times O_3)^\perp\times E$ and $\mathcal{L}_F := (O_2\times O_4)^\perp\times F.\;\;\Box$\\

%\noindent We name the circumcircles of the component triangles and their centers:\\ 
%\centerline{$\mathcal{C}_1$ = circumcircle of $AB\,C$ with center $O_1$, $\mathcal{C}_2$ = circumcircle of $AEF$ with center $O_2$,}\\
%\centerline{$\mathcal{C}_3$ = circumcircle of $CDE$ with center $O_3$, $\mathcal{C}_4$ = circumcircle of $BDF$ with center $O_4$.}
%We first determine the radical center of {$\mathcal{C}_1$, $\mathcal{C}_2$ and $\mathcal{C}_3$. It is the common intersection point of the three radical lines\\
%\centerline{$\mathcal{L}_A := (O_1\times O_2)^\perp\times A = \big(\big((A{-}B)\times(A{-}C)\big)\times \big((A{-}E)\times(A{-}F)\big)\big)\times A$,}
%\centerline{$\mathcal{L}_B := (O_1\times O_3)^\perp\times B = \big(\big((A{-}B)\times(A{-}C)\big)\times \big((A{-}E)\times(A{-}F)\big)\big)\times B$,}
%\centerline{$\mathcal{L}_C := (O_1\times O_4)^\perp\times C = \big(\big((A{-}B)\times(A{-}C)\big)\times \big((A{-}E)\times(A{-}F)\big)\big)\times C$.}
%We calculate\\
%\centerline{  $\mathcal{L}_A \times \mathcal{L}_B = [ef(d(m{-}n){\,-\,}1):df(1{\,-\,}e\,(l{-}n)):de\,(f (l{-}m){\,-\,}1)] = \mathcal{L}_A \times \mathcal{L}_C$.}\\
%The point $M\!q:=[ef(d(m{-}n){\,-\,}1):df(1{\,-\,}e\,(l{-}n)):de\,(f (l{-}m){\,-\,}1)]$
%is also a point on the lines $\mathcal{L}_D := (O_3\times O_4)^\perp\times D$, $\mathcal{L}_E := (O_2\times O_3)^\perp\times E$ and $\mathcal{L}_F := (O_2\times O_4)^\perp\times F.\;\;\Box$\\

According to the principle of duality in projective geometry, the dual version of the theorem is also a valid theorem.\vspace*{1 mm}\\
\noindent\textit{Theorem} 1b.
Each of the four component trilaterals of a complete quadrangle $\{A,B,C,D\}$ has a tangent circle such that any two of these circles intersect at two points and all external similarity points (there are six in total) lie on a line. We call this line the \textit{Miquel-Steiner line} of the quadrangle. See Figure 3.\vspace*{-2mm}\\ 
\begin{figure}[!bhpt]
\includegraphics[height=5.0cm]{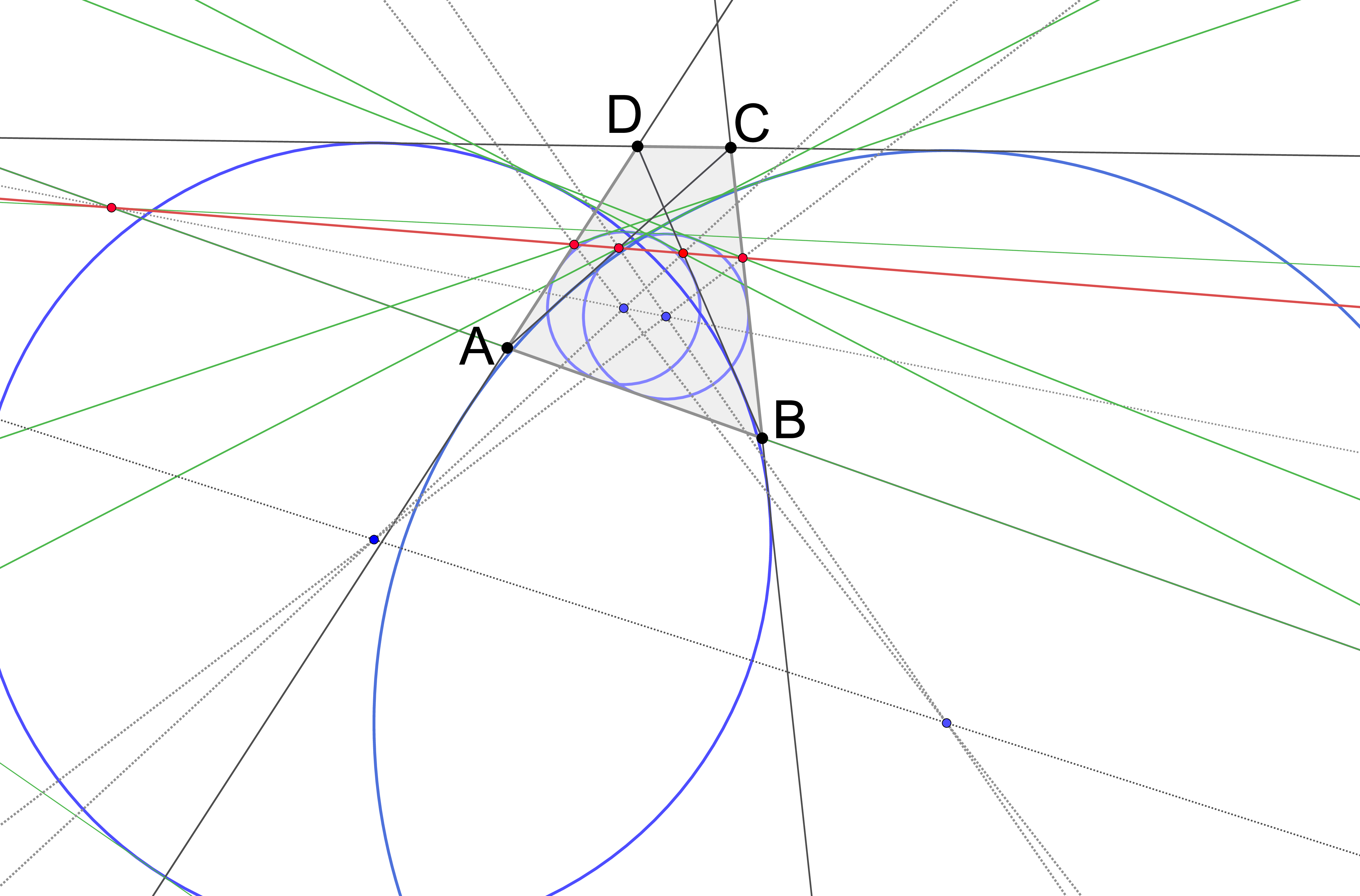}
\centering
\caption{The red line is the Miquel-Steiner line of a quadrangle $\{A,B,C,D\}$ in the euclidean plane. The four blue circles form the only proper quadruple of tangent circles of this quadrangle.} 
%The four pink circles do not; there is no pink circle that touches the line $A\times C$}
\vspace*{-2 mm}\end{figure}

\noindent\textit{The Miquel-Steiner point of a tetragon}.  Let $ABCD$, $A=[1{:}0{:}0], B=[0{:}1{:}0],C=[0{:}0{:}1],$ $D=[d_A{:}d_B{:}d_C]$ (with $d_A > 0)$, be a tetragon. The diagonal points are
$P_1{\,:=\,}(A\times D)\times(B\times C) =[0{:}|d_B|{:}\chi(0,d_B,d_C)\,d_C], P_2{\,:=\,}(B\times D)\times(A\times C)=[d_A{:}0{:}d_C]$, $P_3{\,:=\,}(C\times D)\times(A\times B)=[d_A{:}d_B{:}0].$ We will determine the Miquel-Steiner point of this tetragon, which is defined as the Miquel-Steiner point of the  associated quadrilateral $\boldsymbol{Q\!L}{\,=\,}\{A{\times}B,$ $B{\times}C,C{\times}D,D{\times}A\}$. 
For this purpose, we introduce real numbers\\ 
\hspace*{3 mm}$\alpha{\,:=\,}\mathrm{sgn}(d_B)\sqrt{|d_B^2{+}2r_Ad_B d_C{+}d_C^2|}${\,,\,}$\beta{\,:=\,}\sqrt{|d_A^2{+}2r_Bd_B d_C{+}d_C^2|},$\\
\hspace*{3 mm}$\gamma{\,:=\,}\sqrt{|d_B^2{+}2r_Cd_Bd_C{+}d_C^2|}$\,,\,
{$\delta{\,:=\,}\sqrt{|d_A^2{+}d_B^2{+}d_C^2{+}2r_Ad_Bd_C{+}2r_Bd_Cd_A{+}2r_Cd_Ad_B|}\,.$\vspace*{1 mm} \\ 
The calculation of the Miquel-Steiner point is very similar to that in the proof of Theorem\nobreakspace 1. The result is\\
\centerline{$M\!q_{\scriptscriptstyle{ABCD}}=[d_A \delta\gamma(\alpha(d_B{+}d_C)-1){\,:\,}d_B(\delta(d_B\alpha\gamma{-}\alpha{-}\gamma)+\alpha\gamma){\,:\,}d_C \delta\alpha(\gamma(d_A{+}d_B)-1)]$.}\vspace*{2 mm}\\
We also determine the Miquel-Steiner points of the tetragons $BCDA$ and $ADBC$:\\
\centerline{$M\!q_{\scriptscriptstyle{ABDC}}=[d_A(\delta (d_A\beta\gamma{-}\beta{-}\gamma)+\beta\gamma){\,:\,}d_B \delta\gamma(\beta (d_A{+}d_C)-1){\,:\,}d_C \delta \beta(\gamma(d_A{+}d_B)-1)]$;}\\
\centerline{$M\!q_{\scriptscriptstyle{ADBC}}=[d_A \delta\beta(\alpha(d_B{+}d_C)-1){\,:\,}d_C \delta\alpha(\beta(d_A{+}d_C)-1){\,:\,}d_C(\delta(d_C\alpha\beta{-}\alpha{-}\beta)+\alpha\beta)]$,}\vspace*{0.6 mm}\\
put $M\!q_2:=M\!q_{\scriptscriptstyle{ABCD}}, M\!q_1:=M\!q_{\scriptscriptstyle{ABCD}}, M\!q_3:=M\!q_{\scriptscriptstyle{ADBC}} $,\vspace*{0.6 mm}\\and
call the triple $(M\!q_1{,}M\!q_2{,}M\!q_3)$ the \textit{Miquel-Steiner triangle} of quadrangle $\{A{,}B{,}C{,}D\}$.\vspace*{1 mm}\\

\begin{figure}[!htbp]
\includegraphics[height=6cm]{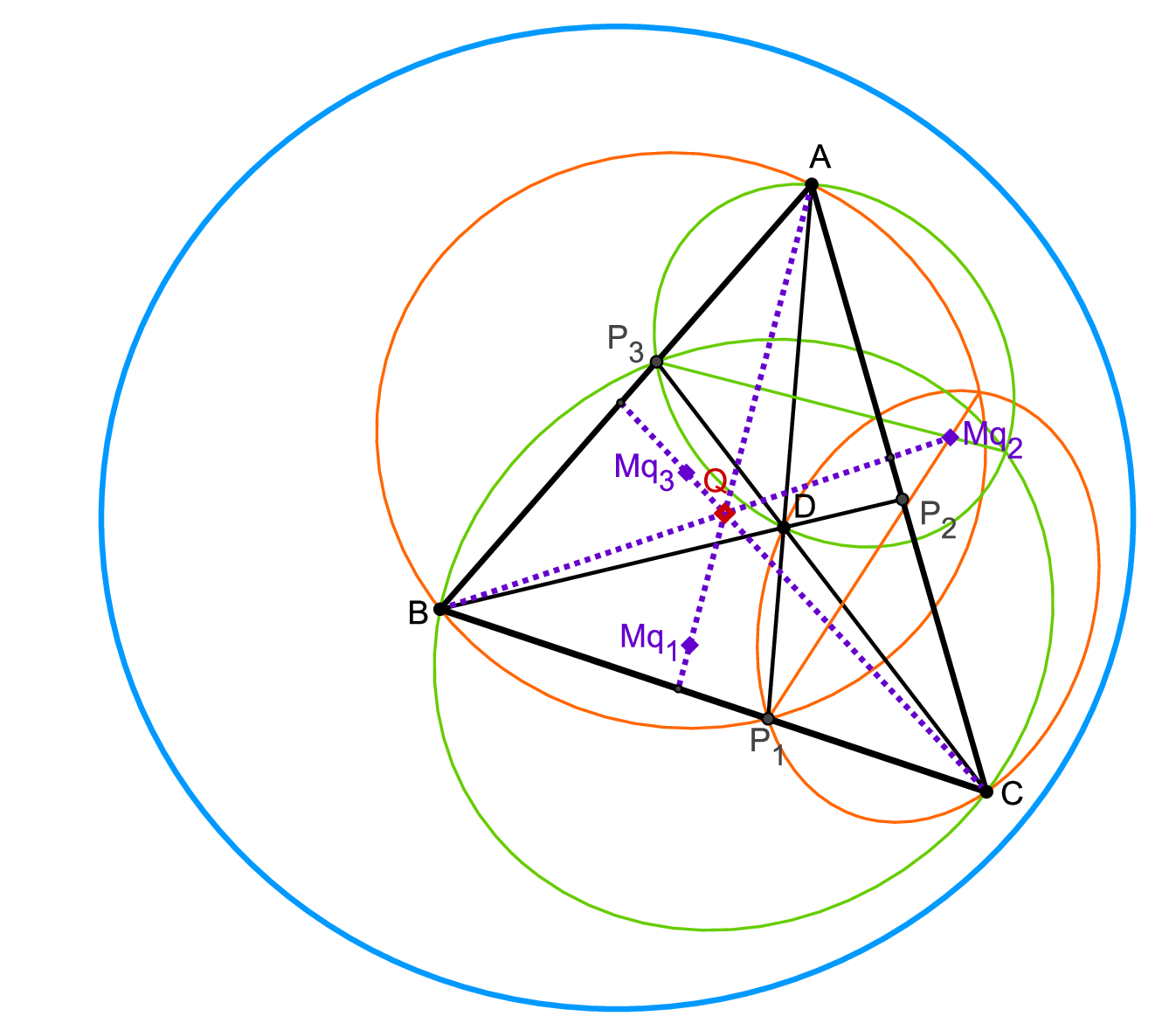}
\centering
\caption{The light blue circle is the absolute conic in a hyperbolic plane.}
\vspace*{0 mm}\end{figure}

\noindent\textit{Odehnal's Miquel-Steiner transformation.} B. Odehnal has shown that in the euclidean plane the Miquel-Steiner triangle of a quadrangle $\{A,B,C,D\}$ is perspective to triangle $ABC$. The mapping that assigns the perspector to point $D$ was called \textit{Miquel-Steiner transformation} by Odehnal and was thoroughly investigated by him, see \cite{Od}.\vspace*{3 mm}\\
\noindent\textit{Theorem} 2. The Miquel-Steiner triangle and triangle $ABC$ are perspectively related to each other also in non-affine planes. 
Using notations introduced above, we obtain the perspector:\\
\centerline{$Q = [d_A \beta\gamma(\alpha(d_B{+}d_C)-1){\,:\,}d_B\alpha\gamma(\beta(d_A{+}d_C)-1){\,:\,}d_C\alpha\beta(\gamma(d_A{+}d_B)-1)]$.}\vspace*{-1 mm}\\

In a euclidean plane, the Miquel-Steiner transformation has exactly one fixed point, and that is the orthocenter of $ABC$.
Experimental constructions using \textit{GeoGebra} show that in non-affine planes the \textit{Miquel-Steiner transformation} has also one fixed point, but it is (in general) not the orthocenter. It can be shown by calculation that the fixed point must satisfy the equations $d_A-d_B =1/{\beta}-1/{\alpha}$ and $d_A-d_C =1/{\gamma}-1/{\alpha}$.\vspace*{2 mm}\\
\noindent\textit{Theorem} 3. Let $ABCD$ be a tetragon with diagonal points $P_1:=(A\times B)\times(C\times D),P_2:= (D\times B)\times(A\times C), P_3:= (B \times C)\times(D\times A)$, and we assume that the Miquel-Steiner point of the tetragon $ABCD$ exists.
\\
Then the Miquel-Steiner point of $ABCD$ is a point on the line $P_1\times P_2$ precisely when the circumcircles of $ABC$ and $CDA$ coincide.\vspace*{2mm}\\
\noindent\textit{Remark}: The corresponding euclidean theorem has been known since the 19th century, with proofs based on Brocard's work.\vspace*{2 mm}\\
\textit{Proof of Theorem} 3.} 
The Miquel-Steiner point of $ABCD$ is \vspace*{0.5 mm}\\
\centerline{$M\!q_{\scriptscriptstyle{ABCD}}=[d_A \delta\gamma(\alpha(d_B{+}d_C)-1){\,:\,}d_B(\delta(d_B\alpha\gamma{-}\alpha{-}\gamma)+\alpha\gamma){\,:\,}d_C \delta\alpha(\gamma(d_A{+}d_B)-1)]$.}\vspace*{1mm}
\noindent This is a point on the line $P_A\times P_C$ iff\vspace*{1 mm}\\  
\hspace*{3mm}$0 = \big(d_A \delta\gamma(\alpha(d_B{+}d_C)-1){\,,\,}d_B(\delta(d_B\alpha\gamma{-}\alpha{-}\gamma)+\alpha\gamma){\,,\,}d_C \delta\alpha(\gamma(d_A{+}d_B)-1)\big)\\
\hspace*{93 mm}\cdot(d_B d_C,-d_A d_C,d_A d_B)$\vspace*{-0.6mm}\\
$\hspace*{5.7mm}= d_A d_B d_C \alpha\gamma (\delta(d_A+d_B+d_C)-1)$.\vspace*{-2 mm}\\

Next, we show that a necessary and sufficient condition for the point $D=[d_A{:}d_B{:}d_C]$ to lie on the circumcircle of  triangle $ABC$ is: $|\delta(d_A+d_B+d_C)|=1$.\vspace*{1mm}\\
\centerline{$D=[d_A{:}d_B{:}d_C]$ is a point on the circumcircle of $ABC$ if and only if}
\hspace*{19.2mm}$0 = 2\,((1{-}r_A)d_Bd_C+(1{-}r_B)d_Cd_A+(1{-}r_C)d_Ad_B)$\\ 
\hspace*{21.0mm}$=(d_A+d_B+d_C)^2 - (r_Ad_Bd_C+r_Bd_Cd_A+r_Cd_Ad_B)$\\ 
\hspace*{21.0mm}$=(d_A+d_B+d_C)^2 - (1/\delta)^2\vspace*{1 mm}$\\
\hspace*{21.0mm}$=(\,d_A+d_B+d_C - 1/\delta\,)(\,d_A+d_B+d_C + 1/\delta\,)$.\vspace*{-2 mm}\\

If $1/\delta=d_A+d_B+d_B$, $B$ is a point on $\Delta_0(CDA)$. However, if $1/\delta=-d_A-d_B-d_B$, $B$ is a point on $\Delta_3(CDA)$, but not on $\Delta_0(CDA)$.$\;\;\;\;\Box$\vspace*{1 mm}\vspace*{1 mm}\vspace*{-3 mm}\\
%\begin{figure}[!hbtp]
%\includegraphics[height=7cm]{Figure4.eps}
%\centering
%\caption{}
%\vspace*{1 mm}\end{figure} \vspace*{2 mm}

\textit{Remark:} In euclidean geometry, the Miquel-Steiner point is obtained by the orthogonal projection of the circumcenter of $ABC$ onto the line $P_1P_3$.
 This does not generally apply to non-affine planes.\vspace*{-2 mm}\\

\noindent\textit{Preliminary remarks} on the next corollary (the dual of Theorem 3).\\
1. A \textit{tangent quadrilateral} is a convex tetragon whose sides are tangent to a single circle.\\
2. The Miquel-Steiner line of a tetragon $ABCD$ is defined as the Miquel-Steiner line of the associated quadrangle $\{A,B,C,D\}$. \vspace*{2 mm}\\
\noindent\textit{Theorem} 3b. 
If a Miquel-Steiner line of a tetragon $ABCD$ exists, then it passes through the diagonal point $P_2=(B\times D)\times(A\times C)$ if and only if $ABCD$ is a tangent quadrilateral and the bisectors of the angles  $\angle(ABC),\angle(CDA)$ meet at the center of the inscribed circle.\vspace*{-2 mm}\\

\section{The Miquel-Steiner theorem in metric-affine planes.}

Following the notation in Figure 1, let the vertices of the complete quadrilateral be labeled accordingly. We let the line at infinity be any line avoiding the vertices. Its tripole with respect to $ ABC$ is the centroid $G$. Now we choose an anisotropic point $K = [k_1{:}k_2{:}k_3]$ (where $k_1k_2k_3 \neq 0$) to be the Lemoine point of $ABC$. In other words, we select the circumconic $k_1 x_2 x_3 + k_2 x_1 x_3 + k_3 x_1 x_2 = 0$ to serve as the circumcircle of $ABC$. 
This fixes the points of intersection of all circles with the line $\mathcal{L}_\infty\!$\,, but also defines the geometry on this line, and even the geometry of the entire plane. If there are two distinct real points of intersection, the plane is a Minkow\-ski plane. If there is only one real intersection point, the plane is galilean. And if there are no real intersection points, the geometry is euclidean.\vspace*{2 mm} \\
We will proof the following (slightly modified version) of the\vspace*{1 mm}\\
\textit{Miquel-Steiner Theorem}: The circumcircles of the four component triangles either touch each other at a real circular point of the plane or they all intersect transversally at an anisotropic point. In any case, their common point is called the \textit{Miquel-Steiner point} of the quadrilateral.\vspace*{2 mm} \\
\hspace*{2 mm}We give two different proofs for this theorem, but since in a galilean plane all circles touch the line at infinity at the absolute pole, we only consider euclidean and Minkow\-ski planes.\vspace*{1 mm}\\
We base the first proof on the following \vspace*{2 mm}\\
\textit{Proposition / Theorem}:
For any circle $\mathcal{C}$ in the plane with a center other than $O$, there exists a triple of real numbers $(p,q,r)$ such that\\ \centerline{$\mathcal{C}= \{[x{:}y{:}z] | k_1yz+k_2zx+k_3xy+(px+qy+rz)(x+y+z)=0\}$.}\\
 The line with the equation $px+qy+rz=0$ is the radical axis of $\mathcal{C}$ and the circumcircle. \vspace*{1 mm} \\
See \cite{Y} for a proof of this theorem.\vspace*{2 mm} \\
\textsl{First proof} of the Miquel-Steiner theorem.\\
We name the vertices of the complete quadrilateral as shown in Figure 1, and the circumcircles of the component triangles as follows:\\ 
\centerline{$\mathcal{C}_{\;}$ = circumcircle of $AB\,C$ , $\mathcal{C}_1$ = circumcircle of $AEF$,}\\
\centerline{$\mathcal{C}_2$ = circumcircle of $BDF$ , $\mathcal{C}_3$ = circumcircle of $CDE$.}\\
We can easily determine the equations for the circles $\mathcal{C}_i, i=1,2,3$, using the proposition above:\\
\centerline{$\mathcal{C}_1=\{[x{:}y{:}z]\} | k_1 y z + k_2 x z + k_3 x y + (x + y + z) (0 x + \dfrac{m k_3}{l - m} y +  \dfrac{n k_2}{l - n} z) = 0$,\;}\\
\centerline{$\mathcal{C}_2=\{[x{:}y{:}z]\} | k_1 y z + k_2 x z + k_3 x y + (x + y + z) (\dfrac{l k_3}{m - l} x + 0 y +  \dfrac{n k_1}{m - n} z) = 0$,}\\
\centerline{$\mathcal{C}_3=\{([x{:}y{:}z]\} | k_1 y z + k_2 x z + k_3 x y + (x + y + z) (\dfrac{l k_2}{n - l} x +  \dfrac{m k_1}{n - m} y + 0 z) = 0$.}\vspace*{1 mm}\\
Now it can be easily checked that the point $M\!q:=(\dfrac{mn}{m-n} k_1 :\dfrac{nl}{n-l}k_2:\dfrac{lm}{l-m}k_3)$ lies on all four circles.\\
In a Minkow\-ski plane, but not in a euclidean plane, the point $M\!q$ may lie on the line at infinity.\\
The statement of the theorem regarding the intersection multiplicities of the circles at $M\!q$ follows from Bezout's theorem. $\;\;\;\Box$\\ 

\hspace*{2 mm} We give a second proof of the Miquel-Steiner theorem by transferring the proof published by Jean-Pierre Ehrmann for the euclidean plane (see \cite{Ehr}) into a proof for the Minkow\-ski plane. To achieve this, we first need to find a theorem for the Minkow\-ski plane that is adequate for the inscribed circle theorem of the euclidean plane.\vspace*{1 mm}\\
\textit{The inscribed theorem for the euclidean space}. 
Let $P,Q,R_1,R_2$ be four distinct points on a circle $\mathcal{C}$. If $R_1, R_2$ are in the same connected component of
$\mathcal{C}\,\mysetminus\,\{P,Q\}$, then $\angle(PR_1Q) = \angle(PR_2Q)$; otherwise, $\angle(PR_1Q) + \angle(PR_2Q) = \pi$. \newpage
 What is meant by the angle measure $\angle $ here? In euclidean geometry, the measure of an angle corresponds to the length of the line segment that the angle cuts out of the line at infinity. This line is an elliptic line; therefore, the length is a real number in the range from $0$ to $\pi$.\\
\hspace*{2mm} In Minkow\-ski geometry, we also define the measure of an angle as the length of the line segment that the angle cuts out of the line at infinity. In this case, the line is a hyperbolic line. But if we adopt a length measurement of segments of a hyperbolic line as described in \cite{Ev1,Ev2}, the inscribed angle theorem in this form can be transferred directly into Minkow\-ski geometry. \vspace*{1 mm}\\
\hspace*{2mm} We should also mention that in Minkow\-ski geometry (as in euclidean geometry) the measures of the interior angles of a triangle sum up to $\pi$.\vspace*{1 mm}\\
\begin{figure}[!thbp]
\includegraphics[height=5.0cm]{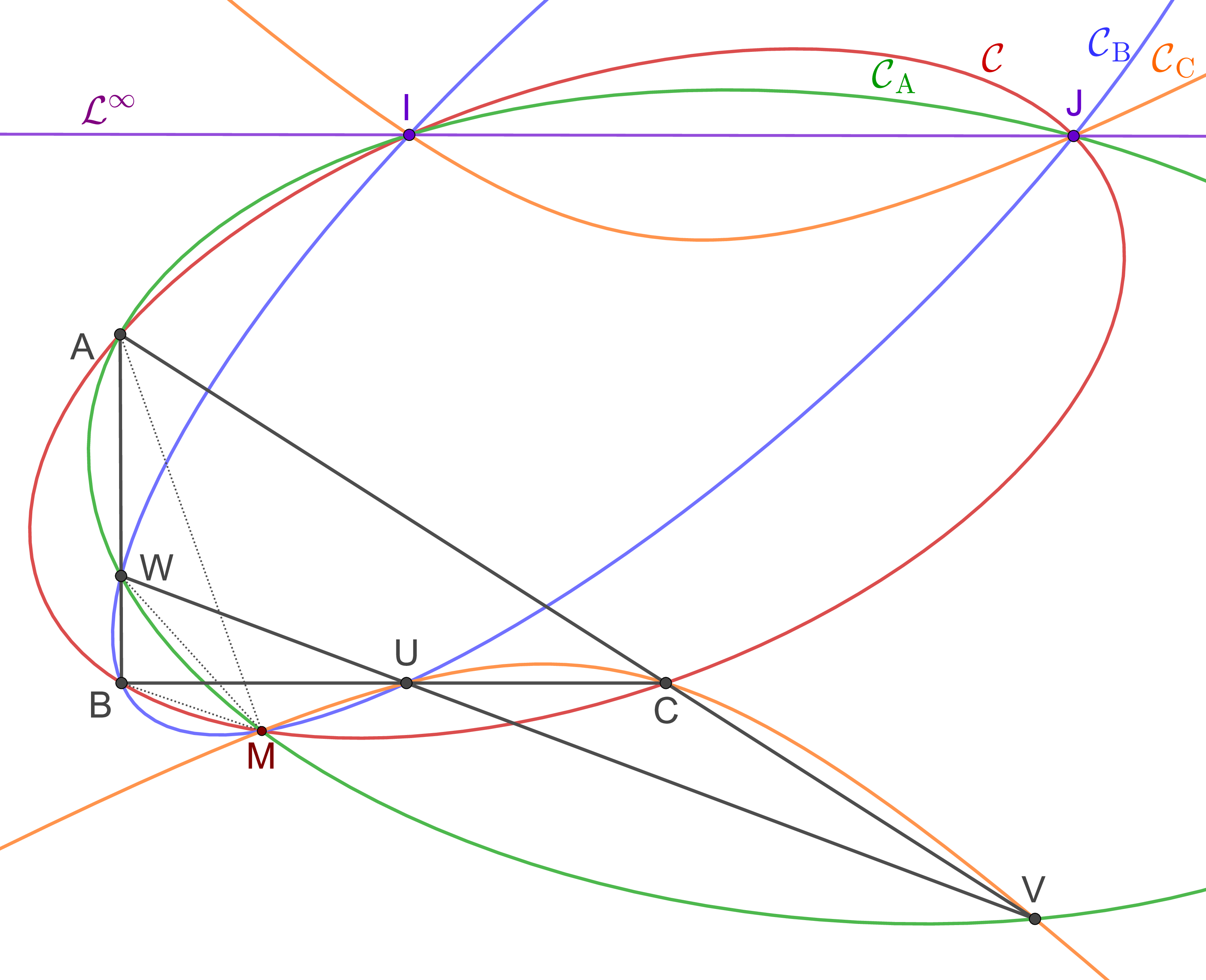}
\centering
\caption{}
\vspace*{-2 mm}\end{figure} 
\noindent\hspace*{-1mm}\textit{Second proof} of the Miquel-Steiner theorem (cf. \cite{Ehr}):\vspace*{0 mm} \\
Points, circles and angles are shown in Figure 5.
Assume that the circles $\mathcal{C}$ and $\mathcal{C}_a$ meet at two distinct points in the affine part of the plane.
Let $M$ be the intersection point of circles $\mathcal{C}$ and $\mathcal{C}_a$, different from $A$.
Because of\\
\hspace*{15.4 mm}$\measuredangle(WMB) = \measuredangle(AMB)+\measuredangle(WMA)$  \\
\hspace*{32 mm}$=\measuredangle(BCA)+\measuredangle(AVW) $ \\                                          
\hspace*{32 mm}$=\measuredangle(ACB) + \measuredangle(WVA)$\\
\hspace*{32 mm}$=\pi-\measuredangle(BCV) - \measuredangle(CVU)$\\
\hspace*{32 mm}$= \measuredangle(VUC)$ \\
\hspace*{32 mm}$= \measuredangle(WUB)$,\\
$M$ is a point on $\mathcal{C}_b$. Accordingly, it can be shown that $M$ lies on $\mathcal{C}_c$.\;\; 
$\Box$\vspace*{2.5mm} \\
\hspace*{2 mm} In 1828, J. Steiner published a short paper \cite{St} with ten problems (questions propos\'{e}es) on a complete quadrilateral. Several proofs of Steiner's conjectures have been published until today, cf. \cite{Cl, Ehr, Fl}. The first proof of the so-called Miquel-Steiner theorem was published by Miquel in 1838, cf. \cite{Mi}. Apart from this one theorem, the other nine problems presented by Steiner (for euclidean geometry) cannot be directly transferred to regular Cayley-Klein planes.\vspace*{1 mm}\\
\hspace*{2 mm} We will now prove a lemma whose validity for the euclidean plane is well-known, but which is not included in Steiner's list. An extension of this lemma to the non-affine planes is not possible.\vspace*{2.5 mm} \\ \newpage
\noindent\textit{Lemma} (cf. \cite{CG,Ro}): For a complete quadrilateral with anisotropic Miquel-Steiner point, as shown in Figure 6, the angles $\angle(C M\!q F)$, $\angle(B M\!q E)$, $\angle(D M\!q A)$ share the same angle bisectors.\vspace*{-2 mm}\\

\begin{figure}[!htbp]
\includegraphics[height=5.5cm]{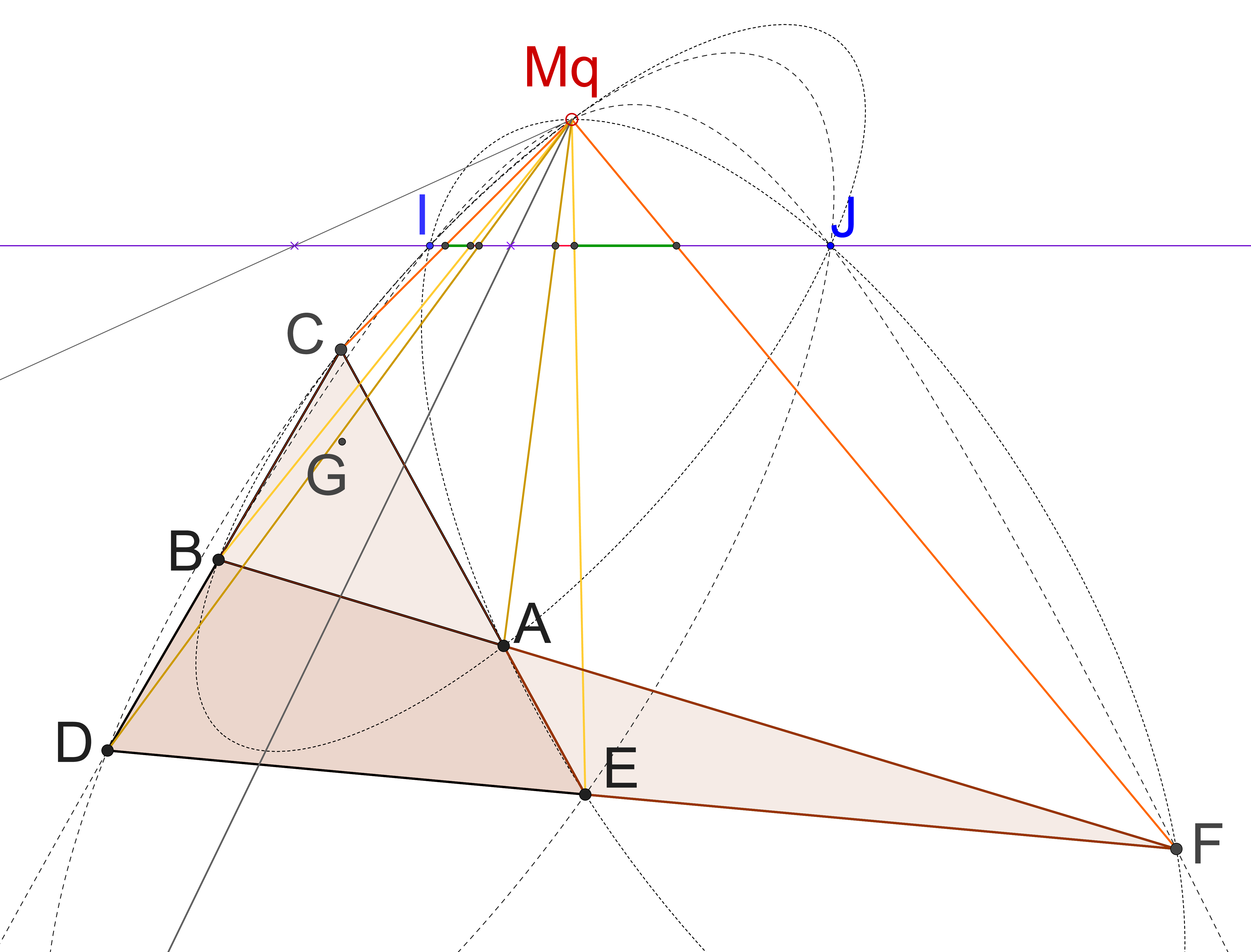}
\centering
\caption{}
\end{figure} \vspace*{-2 mm}

\noindent\textit{Proof}: We choose coordinates of the points $A,B,C,G,D,E,F,K,M\!q$ as given in the first proof of the Miquel-Steiner theorem, and calculate the two circular points $I$ and $J$ on the line at infinity:\\
\centerline{$I = \mathbb{R}\boldsymbol{i}, \boldsymbol{i}=(\rho-k_1{+}k_2{-}k_3,-\rho+k_1{-}k_2{-}k_3,2k_3)$,}\\ 
\centerline{$J = \mathbb{R}\boldsymbol{j}, \boldsymbol{j}=(-\rho-k_1{+}k_2{-}k_3,\rho+k_1{-}k_2{-}k_3,2k_3)$,}
\centerline{with $\rho = \sqrt{k_1^2+k_2^2+k_3^2-2k_1k_2-2k_1k_3-2k_2k_3}$.}\\
$\rho$ is a positive real number if $\mathcal{L}_{\infty}$ is a hyperbolic line, and imaginary if $\mathcal{L}_{\infty}$ is elliptic.\vspace*{0.5 mm}\\
\textit{Note}: This representation of points $I$ and $J$ assumes that $k_3$ is not zero. However, if $k_3 = 0$, then at least one of the two numbers $k_1, k_2$ is not zero, and a corresponding representation of points $I$ and $J$ can be found.\vspace*{0.8 mm}\\
Put $\tilde{A}:= (A\times M\!q)\times \mathcal{L}_{\infty}$, $\tilde{B}:= (B\times M\!q)\times \mathcal{L}_{\infty}$,...,
$\tilde{F}:= (F\times M\!q)\times \mathcal{L}_{\infty}$.\vspace*{0 mm}\\
These points have coordinates  $\tilde{A}=\mathbb{R}\boldsymbol{\tilde{a}},\; \boldsymbol{\tilde{a}}=( l(n-m),n(m-l),m(l-n) )$,\\
\hspace*{3mm}$\tilde{B}=\mathbb{R}\tilde{\boldsymbol{b}},\; \tilde{\boldsymbol{b}}=( n(m-l),m(l-n),l(n-m))\;,\;$ $\tilde{C}=\mathbb{R}\tilde{\boldsymbol{c}},\; \tilde{\boldsymbol{c}}=(m(n-l),l(m-n),n(l-m) )$,\vspace*{0.5mm}\\
\hspace*{3mm}$\tilde{D}=\mathbb{R}\boldsymbol{\tilde{d}},\; \boldsymbol{\tilde{d}}=( (l-m)(l-n)(m-n),l^2n+lm(m-3n)+mn^2,-l^2m+ln(3m-n)-m^2n )$,\vspace*{0.5mm}\\
\hspace*{3mm}$\tilde{E}=\mathbb{R}\boldsymbol{\tilde{e}},\; \boldsymbol{\tilde{e}}=( -\,l^2m+ln(3m-n)-m^2n,(l-m)(l-n)(m-n),l^2n+lm(m-3n)+mn^2 )$,\vspace*{0.5mm}\\
\hspace*{3mm}$\tilde{F}=\mathbb{R}\boldsymbol{\tilde{f}},\; \boldsymbol{\tilde{f}}=( -l^2n+lm(3n-m)-mn^2,l^2m+ln(n-3m)+m^2n,(l-m)(l-n)(n-m) )$.\vspace*{1mm}\\
We will show that $(I,J;\tilde{B},\tilde{D};\tilde{A},\tilde{E})$ and $(I,J;\tilde{C},\tilde{D},\tilde{A},\tilde{F})$ are quadrilateral sets.\vspace*{1mm}\\
We define a function $\tau:(\mathbb{C}^3)^7\to \mathbb{C}$ by\;\;$\tau(\boldsymbol{q},\boldsymbol{r},\boldsymbol{s},\boldsymbol{t},\boldsymbol{u},\boldsymbol{v},\boldsymbol{w}) {:=}$\\ 
%$\hspace*{20mm}\tau(\boldsymbol{q},\boldsymbol{r},\boldsymbol{s},\boldsymbol{t},\boldsymbol{u},\boldsymbol{v},\boldsymbol{w}):=\hspace*{4mm}
%\det(\boldsymbol{q},\boldsymbol{t},\boldsymbol{u})\,\det(\boldsymbol{q},\boldsymbol{v},\boldsymbol{s})\,\det(\boldsymbol{q},\boldsymbol{w},\boldsymbol{r})$\\
%\hspace*{52mm}$-\;\det(\boldsymbol{q},\boldsymbol{v},\boldsymbol{w})\,\det(\boldsymbol{q},\boldsymbol{t},\boldsymbol{s})\,\det(\boldsymbol{q},\boldsymbol{u},\boldsymbol{r})$.\vspace*{1 mm} \\ 
\centerline{$\det(\boldsymbol{q},\boldsymbol{t},\boldsymbol{u})\,\det(\boldsymbol{q},\boldsymbol{v},\boldsymbol{s})\,\det(\boldsymbol{q},\boldsymbol{w},\boldsymbol{r})-\det(\boldsymbol{q},\boldsymbol{v},\boldsymbol{w})\,\det(\boldsymbol{q},\boldsymbol{t},\boldsymbol{s})\,\det(\boldsymbol{q},\boldsymbol{u},\boldsymbol{r})$.}
The statement of the theorem is equivalent to (cf. \cite{Li, RG})\vspace*{1mm}\\
\centerline{$\tau(\boldsymbol{q},\boldsymbol{i},\boldsymbol{j},\tilde{\boldsymbol{b}},\tilde{\boldsymbol{d}},\tilde{\boldsymbol{a}},\tilde{\boldsymbol{e}})=0$ and $\tau(\boldsymbol{q},\boldsymbol{i},\boldsymbol{j},\tilde{\boldsymbol{c}},\tilde{\boldsymbol{d}},\tilde{\boldsymbol{a}},\tilde{\boldsymbol{f}})=0 ,$}\vspace*{1mm}
where $\boldsymbol{q}=(q_1,q_2,q_3)$ can be vector $\boldsymbol{a}$ or any other vector with $q_1{+}q_2{+}q_3 \ne 0$. 
A calculation of $\tau(\boldsymbol{q},\boldsymbol{i},\boldsymbol{j},\tilde{\boldsymbol{b}},\tilde{\boldsymbol{d}},\tilde{\boldsymbol{a}},\tilde{\boldsymbol{e}})$  and of $\tau(\boldsymbol{q},\boldsymbol{i},\boldsymbol{j},\tilde{\boldsymbol{c}},\tilde{\boldsymbol{d}},\tilde{\boldsymbol{a}},\tilde{\boldsymbol{f}})$ shows that both expressions have the common factor $k_1mn(l-m)(l-n)+k_2ln(l-m)(n-m)+k_3lm(l-n)(m-n)$. But $k_1mn(l-m)(l-n)+k_2ln(l-m)(n-m)+k_3lm(l-n)(m-n)=0$, since $M\!q$ is a point on the circumcircle of triangle $ABC$. $\;\;\Box$\vspace*{3 mm}

\section{The Miquel-Steiner theorem for singular circles}

In the preceding sections, circles were assumed to be regular conics. However, it is entirely legitimate to consider singular circles as well, see \cite{Ev0,Ev2,Li}. Since our focus is on the Miquel-Steiner theorem, we are interested here only in circles passing through the vertices of triangles in \(\mathrm{P}\mathbf{U}\). Without attempting a general proof of the Steiner-Miquel theorem for non-regular circles, we demonstrate its validity using three examples.\vspace*{-1mm}\\

Example 1: 
An affine triangle in a Minkowski plane has a singular circumcircle if at least one sideline is isotropic, passing through a circular point. Here, the circumcircle consists of that isotropic line and a second line connecting the other circular point to the opposite vertex. See Figure 7.\vspace*{-1 mm}\\

Example 2:
Let $ABC$ be an affine triangle. By adding the tripolar of its centroid $G$, the affine plane is projectively completed. We now consider the triangle $\{[p_A{:}p_B{:}p_C]{\,\in\,}\mathrm{P}\mathbf{U}\,|$ $\,p_A{\geq}\,0\;\mathrm{and}\;p_B p_C \leq 0\}$. We define the (singular) circumcircle of this triangle as the line $B\times C$ together with the line parallel to it passing through A. See Figure 8.\vspace*{-1 mm}\\

Example 3:
In a hyperbolic plane, a singular circle passing through the three anisotropic vertices of a triangle always consists of two distinct tangents to the absolute conic. Their intersection point is the center of the circle.
See Figure 9.\vspace*{1 mm}\\

\begin{figure}[!hptb]
\includegraphics[height=8.5cm]{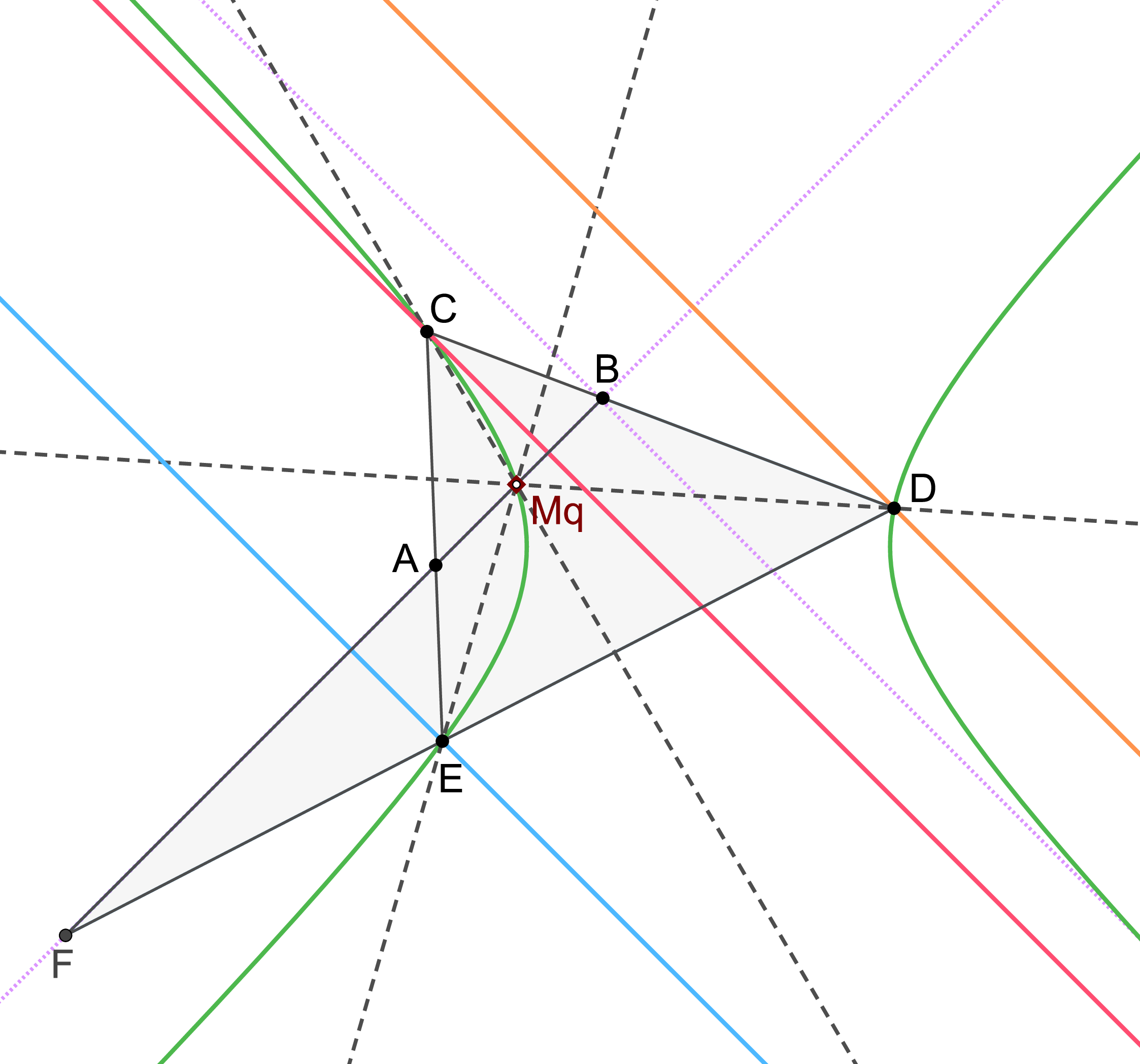}
\centering
\caption{The Miquel-Steiner point of a complete quadrilateral in a Minkowski plane. The line through $A,B,F$ is isotropic.}
\vspace*{-1 mm}\end{figure} \hspace*{-2 mm}

\begin{figure}[!phbt]
\includegraphics[height=8cm]{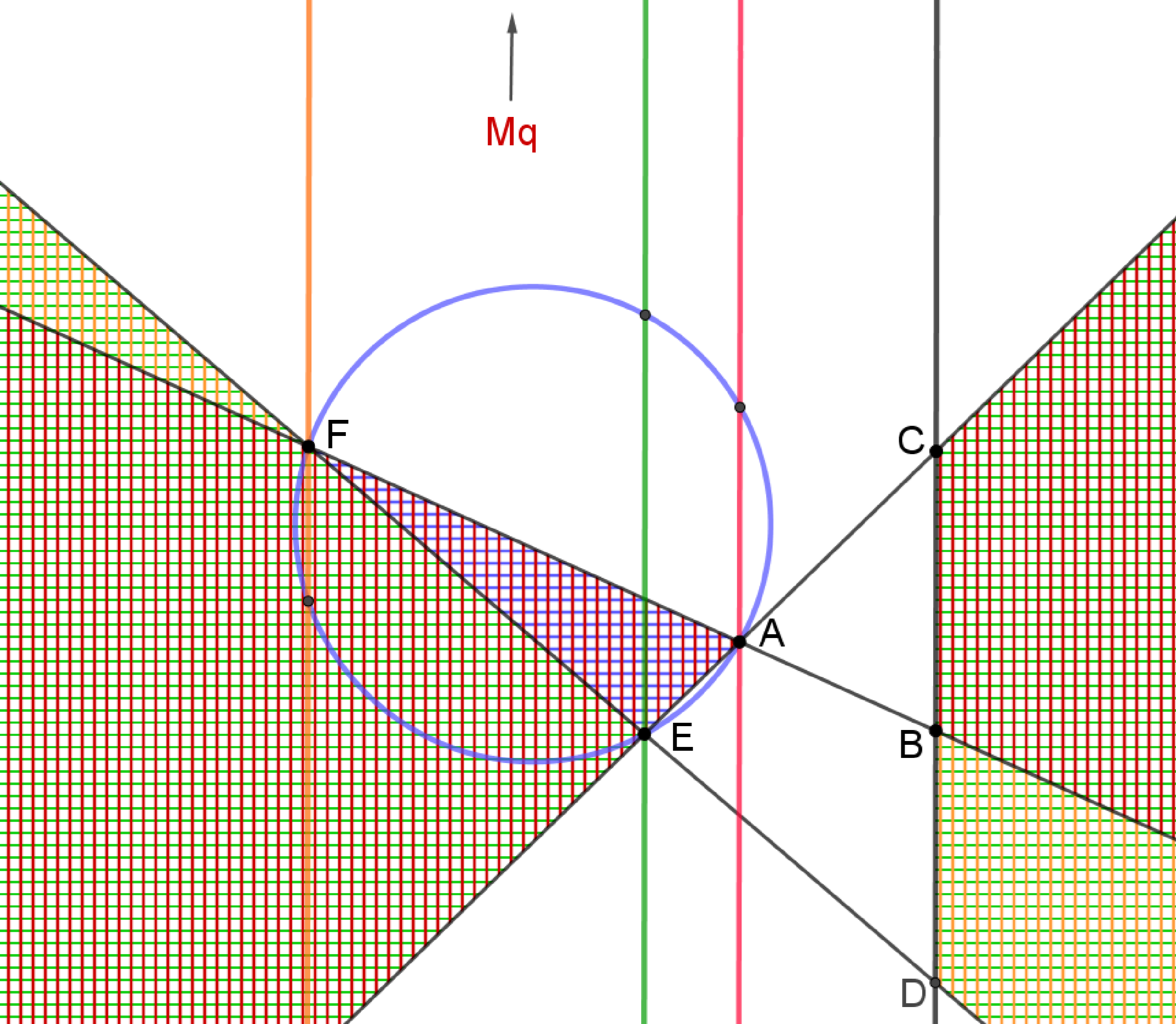}
\centering
\caption{The Miquel-Steiner point of a complete quadrilateral in a euclidean plane. The triangles belonging to the quadrilateral are hatched.}
\vspace*{0 mm}\end{figure} \hspace*{0 mm} 

\begin{figure}[!hptb]
\includegraphics[height=8.0cm]{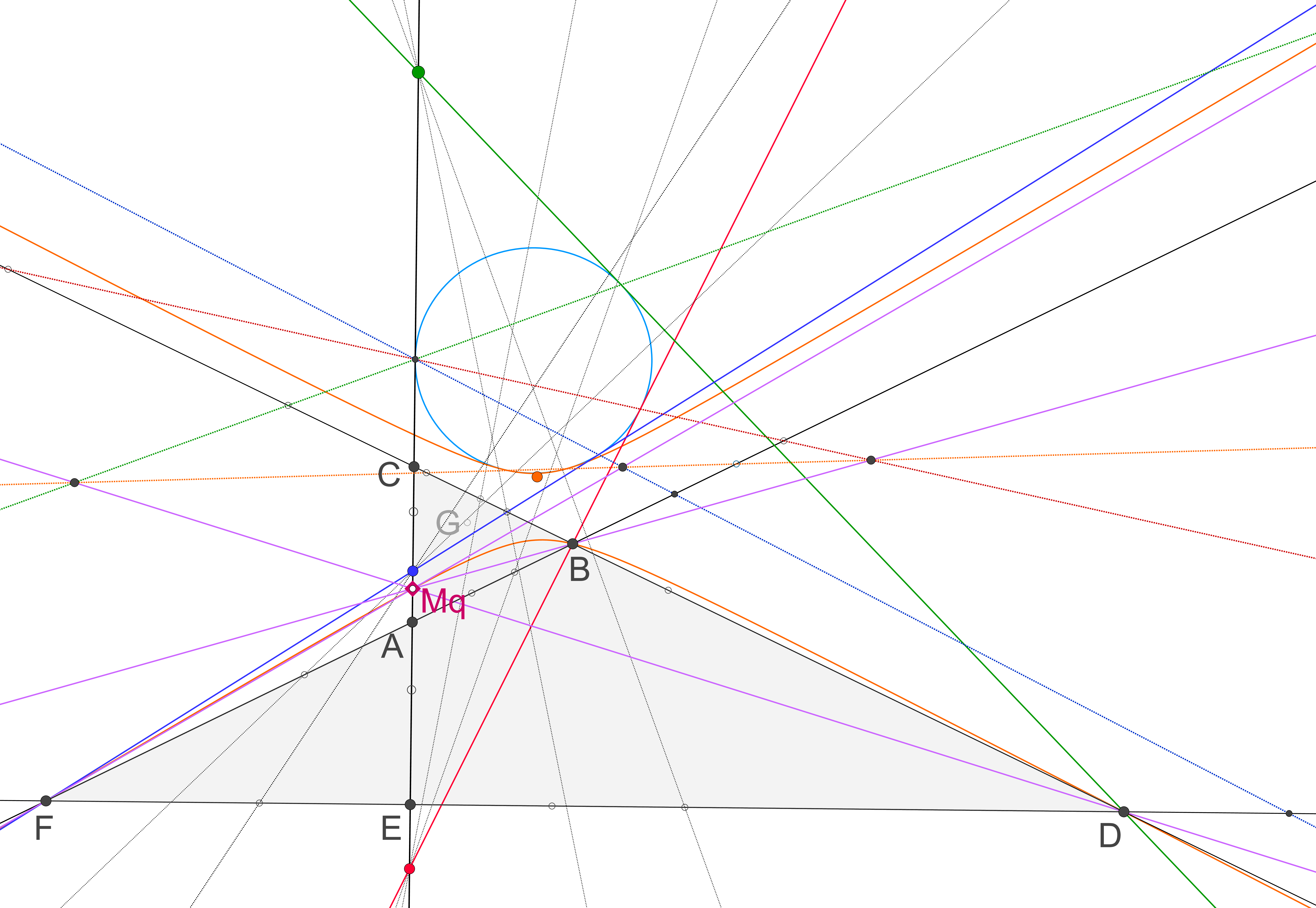}
\centering
\caption{The Miquel-Steiner point of a complete quadrilateral is shown in a hyperbolic plane with the absolute conic in light blue. The line through $C,A,E$ is isotropic. In addition to a regular circle (orange), three of the four circles belonging to the quadrilateral are singular. Their centers (red, blue, green) lie on the line $A \times C$.}
\vspace*{-1 mm}\end{figure} \hspace*{-2 mm}
\newpage

\noindent\textbf{Problem}:\\
In a Cayley-Klein space of dimension $n \geq 2$,  let $\mathcal{H}$ be a set of $n{+}2$ anisotropic hyperplanes in general position. Every subset of $n{+}1$ hyperplanes from $\mathcal{H}$ bounds a simplex, yielding a collection of $n{+}2$ simplices. Assuming each of these simplices possesses a circumsphere, prove or disprove that these $n{+}2$ circumspheres share a common radical center.\\


\begin{thebibliography}{99}
\bibitem {Cl} J.W. Clawson,The Complete Quadrilateral, Ann. Math., 2,20 (1919), 232-261.
\bibitem {CG} H.S.M. Coxeter, S.L. Greitzer, \textit{Geometry Revisited}, The Mathematical Association of America, 1967.
\bibitem {Ehr} J.-P. Ehrmann, Steiner's Theorems on the Complete Quadrilateral, Forum Geom., Vol. 4 (2004), 35-52.
\bibitem {Ev0} M.\,Evers,\;Symbolic Substitution Has a Geometric Meaning, Forum Geom., Vol. 14 (2014), 217-232.
\bibitem {Ev1} M.\,Evers,\;On the Geometry of a Triangle in the Elliptic and in the Extended Hyperbolic Plane,\\ 
\hspace*{1mm}arXiv:1908.11134, 2019.
\bibitem {Ev2} M.\,Evers,\;Geometry on real projective Cayley-Klein spaces, arXiv:2301.04024, 2023.
\bibitem {Ev3} M.\,Evers,\;Quadri-Figures in Cayley-Klein Planes: All Around the Newton Line, arXiv:2409.17802, 2024.
\bibitem {Fl} J. A. Flos,\;Museization of Steiner's 10 Questions on the complete quadrilateral, 2019, available at\\
\hspace*{1mm}https://www.semanticscholar.org
\bibitem {Li} St. Liebscher,\;\textit{Projektive Geometrie der Ebene}, Springer, 2017.
\bibitem {Mi} A. Miquel,\;Théorèmes de Géométrie, Journal de mathématiques pures et appliquées 1re série, tome 3 (1838), p. 485-487.
\bibitem {Od} B. Odehnal,\;A Miquel-Steiner Transformation, KoG, Vol. 27, 2023.
\bibitem {RG} J. Richter-Gebert,\;\textit{Perspectives on Projective Geometry}, Springer, Berlin, 2011.
\bibitem {Ro} V. Rong,\;Complete Quadrilaterals and the Miquel Point, 2021, available at\\
\hspace*{1mm} https://www.lessvrong.com/math/handouts/miquel-summer2021.pdf
\bibitem {Sch} E. Schmidt,\;Miquel-,Poncelet- und Bennett-Punkt eines Vierecks, 2011, available at\\
\hspace*{1mm} http://eckartschmidt.de/Pktve.pdf
\bibitem {St} J. Steiner,\;Questions propos\'{e}es. Th\'{e}or\`{e}me sur le quadrilat\`{e}re complet, tome 18 (1827-1828), p. 302-304.
\bibitem {Th} P. Thurnheer,\;Int. J. Geom., Vol. 11 (2022), No. 2, 67-77.
\bibitem {WO} G. Weiss, B. Odehnal, Miquel’s Theorem and its Elementary Geometric Relatives, KoG, Vol. 28, 2024. 
\bibitem {Y} P. Yiu,\;\textit{Introduction to the Geometry of the Triangle}, Florida Atlantic University Lecture Notes, 2001.		
\bibitem {Z} Y. Zakharyan, Miquel-Steiner's point locus, arXiv:2002.12777v1, 2020.
\bibitem {GG} GeoGebra,\;Ein Softwaresystem f\"ur dynamische Geometry und Algebra, invented by M. Hohenwarter. 
				
\end{thebibliography}
\end{document}